\documentclass[12pt]{amsart}
\usepackage[margin = 1in]{geometry}
\usepackage{amsmath,amssymb,amsthm,mathtools}
\usepackage{colordvi}
\usepackage[dvipsnames]{xcolor}
\usepackage{ulem}  
\usepackage{graphicx}
\usepackage{verbatim}
\usepackage{bbm}

\usepackage{framed}

\newtheorem{theorem}{Theorem}

\newtheorem{corollary}[theorem]{Corollary}

\theoremstyle{definition}

\newcommand{\CC}{\mathbb{C}}
\newcommand{\RR}{\mathbb{R}}
\newcommand{\ZZ}{\mathbb{Z}}

\newcommand{\KK}{\mathbb{K}}

\newcommand{\calV}{\mathcal{V}}

\DeclareMathOperator{\var}{{\rm var}}
\DeclareMathOperator{\Syl}{{\rm Syl}}

\definecolor{TAMU}{RGB}{140,0,0}
\definecolor{myblue}{RGB}{0,0,198}
\definecolor{myred}{RGB}{182,0,0}

\newcommand{\defcolor}[1]{{\color{TAMU}#1}}
\newcommand{\demph}[1]{\defcolor{{\sl #1}}}

\title{Real solutions to systems of polynomial equations in Macaulay2}
\author[J.~Lopez Garcia]{Jordy Lopez Garcia} 
\address{J.~Lopez Garcia\\ 
         Department of Mathematics\\ 
         Texas A\&M University\\ 
         College Station\\ 
         Texas \ 77843\\ 
         USA} 
\email{jordy.lopez@tamu.edu}
\urladdr{https://jordylopez27.github.io/}
\author[K.~Maluccio]{Kelly Maluccio} 
\address{K.~Maluccio\\ 
         Department of Mathematics\\ 
         Austin Community College\\ 
         Austin\\ 
         Texas \ 78752\\ 
         USA} 
\email{kmaluccio@gmail.com}
\urladdr{https://github.com/kmaluccio}
\author[F.~Sottile]{Frank Sottile} 
\address{F.~Sottile\\ 
         Department of Mathematics\\ 
         Texas A\&M University\\ 
         College Station\\ 
         Texas \ 77843\\ 
         USA} 
\email{sottile@tamu.edu} 
\urladdr{https://www.math.tamu.edu/\~{}sottile} 
\author[T.~Yahl]{Thomas Yahl} 
\address{T.~Yahl\\ 
         Department of Mathematics\\ 
         Texas A\&M University\\ 
         College Station\\ 
         Texas \ 77843\\ 
         USA} 
\email{thomasjyahl@tamu.edu} 
\urladdr{https://tjyahl.github.io/} 
\thanks{Research supported in part by Simons Collaboration Grant for Mathematicians 636314.}
\thanks{\texttt{RealRoots} version 0.1}

\begin{document}

\begin{abstract}
 The \textit{Macaulay2} package \texttt{RealRoots} provides symbolic methods to study real solutions to systems of polynomial equations.
 It updates and expands an earlier package developed by Grayson and Sottile in 1999.
 We provide mathematical background and descriptions of the  \texttt{RealRoots} package, giving examples which illustrate some of its
 implemented methods.
 We also prove a general version of Sylvester's Theorem whose statement and proof we could not find in the literature.
\end{abstract}

\maketitle

\section*{Introduction}

Understanding the number of real solutions to systems of polynomial equations is fundamental for real algebraic geometry and for
applications of algebraic geometry.
In 1999, Grayson and Sottile \cite{So_M2} developed the \textit{Macaulay2} package \texttt{realroots} for this purpose.
That package had limited functionality, was not documented, and not all of its implemented methods remain compatible with
modern releases of \textit{Macaulay2}. 

The \textit{Macaulay2} package \texttt{RealRoots} expands and modernizes \texttt{realroots}, superseding it.
\texttt{RealRoots} implements symbolic methods for studying real solutions to polynomial systems.
This note provides some mathematical background and examples of methods from the package.
Its  three sections each describe related methods.

Section 1 describes methods for counting and isolating real roots of univariate polynomials, as well as methods for determining if a
polynomial is Hurwitz-stable.
We give an extension of Sylvester's Theorem that we could not find in the literature and sketch its proof.

Section 2 describes methods involving elimination that reduce a zero-dimensional system of multivariate polynomials to a univariate
polynomial for solving, studying the number of real solutions, or addressing other arithmetic questions, such as Galois groups.

Section 3 describes a further method for studying zero-dimensional systems based on the trace symmetric form.

\section{Real roots of univariate polynomials}\label{S:one}

Let $f\in\RR[x]$ be a polynomial.
It has the form
 \[
   f\ =\ c_kx^{a_k}  + \dotsb + c_{1}x^{a_{1}} + c_0x^{a_0}\,,
 \]
where $a_k> \dotsb > a_1 > a_0 \geq 0$ are integers and for $0\leq i \leq k$, $c_{i}\in\RR$ is nonzero.
Let $\defcolor{\var(c_0,\dotsc,c_k)}\vcentcolon=\#\{1\leq i\leq k\mid c_{i-1}c_i<0\}$ be the number of variations in sign of the
coefficients of $f$.
Descartes' Rule of Signs \cite{So_Book} gives an upper bound for the number of positive real roots of $f$.

\begin{theorem}[Descartes' Rule of Signs]
  The number, $r$,  of positive real roots of $f$, counted with multiplicity, is at most $\var(c_0,\dotsc,c_k)$ and the difference
  $\var(c_0,\dotsc,c_k)-r$ is even.
\end{theorem}

Given any sequence $c=(c_0,\dotsc,c_k)$, the \demph{variation},  \defcolor{$\var(c)$} of $c$ is the
number of variations in sign after removing all zero terms.
\begin{leftbar}
\begin{verbatim}
i1 : loadPackage("RealRoots");

i2 : variations {2, -3, 0, -8, 12, 0, 0, 8, -12, 0}

o2 = 3
\end{verbatim}
\end{leftbar}
For a sequence of polynomials  $\defcolor{f_\bullet}=(f_0,\dotsc,f_k)$ in $\RR[x]$ and $a\in\RR$, \defcolor{$\var(f_\bullet,a)$} is the
variation in the sequence 
$(f_0(a),\dotsc,f_{k}(a))$. 
We extend this to $a\in\{\pm\infty\}$, by taking $f(\infty)$ to be the leading coefficient of $f(x)$ and $f(-\infty)$ to be the leading
coefficient of $f(-x)$.

Given a polynomial $f\in\RR[x]$ of degree $k$, consider its sequence of derivatives,
 \[
   \defcolor{\delta f}\ \vcentcolon= \left(f(x),f'(x),f''(x),\dots,f^{(k)}(x)\right)\,.
 \]
For $a<b$ in $\RR\cup\{\pm \infty\}$, let \defcolor{$r(f,a,b)$} be the number of roots of $f$ in the interval $(a,b\hspace{.05cm}]$, counted
with multiplicity.
Budan and Fourier~\cite[Ch.\ 2]{So_Book} generalized Descartes' Rule.

\begin{theorem}[Budan-Fourier]
  We have that $r(f,a,b)\leq \var(\delta f,a) -\var(\delta f,b)$, and the difference
  $\var(\delta f,a) -\var(\delta f,b)-r(f,a,b)$ is even. 
\end{theorem}

Descartes' Rule is when $a=0$ and $b=\infty$.
Let us consider an example.
\begin{leftbar}
\begin{verbatim}
i3 : R = QQ[x];

i4 : f = x*(2*x - 3)*(x^4 - 2)^2

       10     9     6      5     2
o4 = 2x   - 3x  - 8x  + 12x  + 8x  - 12x

i5 : BudanFourierBound(f, 0, infinity)

o5 = 3

i6 : BudanFourierBound(f, -2, 1)

o6 = 7
\end{verbatim}
\end{leftbar}
\noindent Note that $r(f,0,\infty)=r(f,-2,1)=3$, as the real roots of $f$ are $-\sqrt[4]{2},-\sqrt[4]{2},0,\sqrt[4]{2},\sqrt[4]{2},3/2$.

In contrast to these bounds, 
Sylvester's Theorem determines the actual number of real roots, and more.
The \demph{Sylvester sequence}, \defcolor{$\Syl(f,g)$} of polynomials $f,g\in\RR[x]$ is the sequence
$\left(f_0,f_1,\dotsc,f_k\right)$ of nonzero polynomials, where $\defcolor{f_0}\vcentcolon= f$, $\defcolor{f_1}\vcentcolon= f'\cdot g$,
and for $i\geq 1$, 
  \[
    \defcolor{f_{i+1}}\ \vcentcolon=\ -1\cdot \mbox{remainder}(f_{i-1},f_i)\,,
  \]
the negative remainder term in the division of $f_{i-1}$ by $f_i$.
The last nonzero remainder is $f_k = \gcd(f,f'g)$.
Observe that for each $1\leq i\leq k$ there exists $\defcolor{q_i}\in\RR[x]$ such that
 \begin{equation}\label{Eq:divisionAlgorithm}
    f_{i-1}\ =\ q_i(x)f_i(x)-f_{i+1}(x)\,.
 \end{equation}
The \demph{reduced Sylvester sequence}, $\defcolor{g_\bullet}=(g_0,\dotsc,g_k)$ is obtained by dividing each term in the Sylvester sequence
by $f_k=\gcd(f,f'g)$, so that $g_if_k=f_i$ for each $i$. 
Note that  $g_k=1$ and elements of the reduced Sylvester sequence satisfy~\eqref{Eq:divisionAlgorithm} with $g_j$ replacing $f_j$.

\begin{theorem}[Sylvester]
  \label{Th:Sylvester}
  Let $f,g\in\RR[x]$ and suppose that $g_\bullet$ is the reduced Sylvester sequence of $f$ and $g$.
  For $a<b$ in $\RR\cup\{\pm\infty\}$ we have
  \begin{eqnarray*}
    \qquad\var(g_\bullet,a)-\var(g_\bullet,b)&=&
    \#\{\zeta\in(a,b\hspace{.05cm}]\mid f(\zeta)=0\mbox{ and } g(\zeta)>0\}  \ -\
      \\
    && \#\{\zeta\in[\hspace{.05cm}a,b)\mid f(\zeta)=0\mbox{ and } g(\zeta)<0\}\,.\qquad
  \end{eqnarray*}
\end{theorem}

Observe the different roles that the endpoints $\{a,b\}$ play in this formula.

\begin{proof}
  In~\cite[Thm.\ 2.55]{BPR}, Sylvester's Theorem is stated and proven 
  when $f$ does not vanish at $a$ or at $b$, and it is in terms of the Sylvester sequence $\Syl(f,g)$.
  That proof proceeds by studying $\var(\Syl(f,g), t)$ as $t$ increases from $a$ to $b$, noting that it may only
  change when $t$ passes a root of some element of the Sylvester sequence.
  Since multiplying a sequence by a nonzero number $f_k(t)$ does not change its variation, the proof in~\cite{BPR}
  establishes this refined  version when $f$ does not vanish at $a$ or at $b$.
  We proceed with the general case.

  Let $g_\bullet$ be the reduced Sylvester sequence of $f$ and $g$.
  The variation $\var(g_\bullet,t)$ may only change when $t$ passes a root $\zeta\in[a,b\hspace{.05cm}]$ of some $g_i$ in
  $g_\bullet$. 
  Observe that $\zeta$ cannot be a root of two consecutive elements of $g_\bullet$.
  If it were, then by~\eqref{Eq:divisionAlgorithm} and induction, it is a root of all elements of $g_\bullet$, and thus of
  $g_k=1$, which is a contradiction.
  Suppose that $g_i(\zeta)=0$ for some $i\geq 1$.
  By~\eqref{Eq:divisionAlgorithm} again, $g_{i-1}(x)$ and $g_{i+1}(x)$ have opposite signs for $x$ near $\zeta$ and thus
  $g_{i-1},g_i,g_{i+1}$ do not contribute to any change in $\var(g_\bullet,t)$ for $t$ near $\zeta$.
  This remains true if $\zeta=a$ and $t$ increases from $a$ or if $\zeta=b$ and $t$ approaches $b$.

  We now suppose that $g_0(\zeta)=0$ and thus $g_1(\zeta)\neq 0$.
  Then we have $f(\zeta)=0$.
  Let $m$ be the multiplicity of the root $\zeta$ of $f$ so that $f=(x{-}\zeta)^m h$ with $h(\zeta)\neq 0$.
  If $g(\zeta)=0$, then $(x{-}\zeta)^m$ divides $f'g$ and thus $f_k$, and so $g_0=f/f_k$ does not vanish at $\zeta$.
  Thus $g(\zeta)\neq 0$.

  Notice that $\defcolor{h_0}\vcentcolon=f/(x{-}\zeta)^{m-1}$ and $\defcolor{h_1}\vcentcolon=f'g/(x{-}\zeta)^{m-1}$ have the same signs for
  $x$ near $\zeta$ as do $g_0$ and $g_1$.
  A computation reveals that $h_1=mhg+(x{-}\zeta)h'g$.
 Choose $\epsilon>0$ so that $\zeta$ is the only root of any element in $g_\bullet$ lying in the interval
 $[\zeta-\epsilon,\zeta+\epsilon]$.
 We have:
 \[
 \begin{array}{c|c|l}
   x & h_0(x) & h_1(x)\\\hline
   \zeta-\epsilon & -\epsilon h(\zeta-\epsilon)  &
        mh(\zeta-\epsilon)g(\zeta-\epsilon) - \epsilon h'(\zeta-\epsilon)g(\zeta-\epsilon)  \rule{0pt}{13pt}\\
   \zeta     &     0    &   mh(\zeta)g(\zeta)  \rule{0pt}{13pt}\\
   \zeta+\epsilon & \epsilon h(\zeta-\epsilon)  &
        mh(\zeta+\epsilon)g(\zeta+\epsilon) + \epsilon h'(\zeta+\epsilon)g(\zeta+\epsilon)  \rule{0pt}{13pt}
 \end{array}
 \]

 Suppose that $g(\zeta)>0$.
 Then the sign of $h_1$ on $[\zeta-\epsilon,\zeta+\epsilon]$ is opposite to the sign of $h_0(\zeta-\epsilon)$, but the same as the sign of
 $h_0(\zeta+\epsilon)$.
 Thus the variation $\var(g_\bullet,t)$ decreases by 1 as $t$ passes from $\zeta-\epsilon$ to $\zeta$, but is unchanged as $t$
 passes from $\zeta$ to $\zeta+\epsilon$.

 Suppose that $g(\zeta)<0$.
 Then the sign of $h_1$ on $[\zeta-\epsilon,\zeta+\epsilon]$ is the same as the sign of $h_0(\zeta-\epsilon)$, but opposite to the sign of
 $h_0(\zeta+\epsilon)$.
 Thus the variation $\var(g_\bullet,t)$ is unchanged as $t$ passes from $\zeta-\epsilon$ to $\zeta$, but increases by 1 as $t$ 
 passes from $\zeta$ to $\zeta+\epsilon$.

 Now consider the variation $\var(g_\bullet,t)$ for $t\in[a,b\hspace{.05cm}]$.
 This may only change at a number $\zeta\in[a,b]$ if $f(\zeta)=0$.
 If $g(\zeta)>0$ and $\zeta\neq b$, then it decreases by 1.
 If $g(\zeta)<0$ and $\zeta\neq a$, then it increases by 1.
 It is otherwise unchanged.
 This completes the proof.
 \end{proof}

The \demph{Sturm sequence} of a polynomial $f\in\RR[x]$ is the Sylvester sequence $\Syl(f,1)$.
The \demph{reduced Sturm sequence} of $f$ is the reduced Sylvester sequence of $f$ and $1$.

\begin{corollary}[Sturm's Theorem]
  Let $f\in\RR[x]$ and $a<b$ in $\mathbb{R}\cup\{\pm\infty\}$.
  Let $g_\bullet$ be the reduced Sturm sequence of $f$.
  Then the number of zeros of $f$ in the interval $(a,b\hspace{.05cm}]$ equals  $\var(g_\bullet,a) - \var(g_\bullet,b)$.
\end{corollary}

Using the reduced Sylvester sequence of $f$ and $-1$, we obtain the number of zeros of $f$ in $[a,b)$. Let us continue with the same polynomial $f=x(2x-3)(x^4-2)^2$ as before.
\begin{leftbar}
\begin{verbatim}
i7 : SylvesterCount(f, x^2 - 1, -2, 3)

o7 = 2

i8 : SturmCount(f)

o8 = 4
\end{verbatim}
\end{leftbar}
\noindent Calling {\tt SturmCount(f)} without endpoints $a,b$ returns the total number of real roots of $f$.

\begin{figure}[htb]
  \centering
  \begin{picture}(220,125)

    \put(0,0){\includegraphics{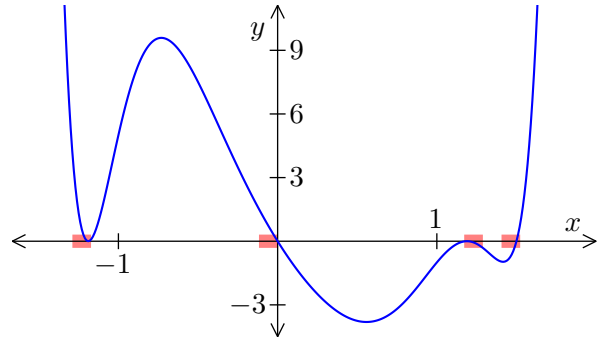}}
    \put(209,40){\small$x$}    \put( 90,114){\small$y$}
    \put( 31,24){\small$-1$}
    \put(157.5,41){\small$1$}  \put( 82, 9){\small$-3$}
    \put(105,57){\small$3$}  \put(105,81){\small$6$} \put(105,105){\small$9$}
        
  \end{picture}
\caption{Graph of $f$.}\label{F:One}
\end{figure}

Figure~\ref{F:One} shows the graph of $f$ in a neighborhood of its real roots.
Note that $x^2-1$ is negative only at the root $0$.

An application of Sturm's Theorem is to give \demph{isolating intervals}, which are disjoint intervals each containing exactly one root
of $f$. 
Our implementation gives a list of pairs $\{p,q\}$ such that $(p,q]$ contains a unique root of $f$ and $q-p$ is less than a
user-provided tolerance.
The numbers $p,q$ are dyadic, lying in $\ZZ[\frac{1}{2}]$, as they are found in a binary search.
\begin{leftbar}
\begin{verbatim}
i9 : realRootIsolation(f, 1/5)

         165    75       15       75  165    45  195
o9 = {{- ---, - --}, {- ---, 0}, {--, ---}, {--, ---}}
         128    64      128       64  128    32  128
\end{verbatim}
\end{leftbar}
\noindent These isolating intervals are shaded in Figure~\ref{F:One}. 

Thomas~\cite{Thomas} observed that recursively iterating Sturm's Theorem on $f_k=\gcd(f,f')$ can be used to give the
number of real roots of $f$, counted with multiplicity.
This same idea may be used to extend Sylvester's Theorem to give the count with multiplicity.
\begin{leftbar}
\begin{verbatim}
i10 : SturmCount(f, -1, 2, Multiplicity => true)

o10 = 4
\end{verbatim}
\end{leftbar}

A polynomial $f\in\RR[x]$ is \demph{Hurwitz-stable} if its complex roots all have negative real parts.
All solutions to the system of constant coefficient ordinary differential equations 
 \[
  \dot{y}\ =\ Ay
 \]
are asymptotically stable ($\lim_{t\to\infty}y(t)=0$) when all eigenvalues $\zeta$ of $A$ have negative real part,
equivalently when the characteristic polynomial of $A$ is Hurwitz-stable.
 
Given a polynomial $f=c_kx^k+c_{k-1}x^{k-1}+\dotsb+c_0$,  let \defcolor{$H$} be the matrix
\[
\left(\begin{matrix}
  c_{k-1} & c_{k-3} & \dotsb & 0  & 0 \\
  c_k    & c_{k-2} &  \dotsb& 0 & 0 \\
  \vdots & \vdots & \ddots&    &\vdots\\
  \vdots & \vdots & & \ddots&\vdots\\
    0    &    0   & \dotsb&c_2 & c_0
\end{matrix}\right)
\]
For $1\leq i\leq k$, the \demph{Hurwitz determinant}, \defcolor{$\Delta_i$} is the $i$th principal minor of $H$.

\begin{theorem}[Hurwitz~\cite{Hurwitz}]
  Suppose that $c_k>0$.
  Then $f$ is Hurwitz-stable if and only if all the Hurwitz determinants $\Delta_{1},\dots,\Delta_{k}$ are positive.
\end{theorem}

\texttt{RealRoots} implements both the Hurwitz matrix and this test for Hurwitz-stability.
\begin{leftbar}
\begin{verbatim}
i11 : HurwitzMatrix(x^4 + 5*x^3 + 7*x^2 + 11*x + 13) 

o11 = | 5 11 0  0  |
      | 1 7  13 0  |
      | 0 5  11 0  |
      | 0 1  7  13 |

               4        4
o11 : Matrix QQ  <--- QQ

i12 : isHurwitzStable(x^4 + 5*x^3 + 7*x^2 + 11*x + 13) 

o12 = false

i13 : isHurwitzStable(x^4 + 9*x^3 + 7*x^2 + 5*x + 3)

o13 = true
\end{verbatim}
\end{leftbar}
%

\section{Elimination}\label{S:two}

Elimination is a classical symbolic method often used to solve systems of equations involving multivariate polynomials.
Geometrically, it gives the image of a variety under a polynomial map, such as a coordinate projection.
\texttt{RealRoots} implements methods for zero-dimension\-al ideals that reduce their study to that of univariate polynomials.

Let  $\KK$ be a field and $I\subseteq\KK[x_1,\dotsc,x_n]$ be a zero-dimension\-al ideal  with scheme $\calV(I)\subseteq\KK^n$.
The Artinian ring $\defcolor{R}\vcentcolon=\KK[x_1,\dotsc,x_n]/I$ is a vector space over $\KK$ of dimension
$\defcolor{d}\vcentcolon=\mbox{degree}(I)$, and $\#\calV(I)\leq d$.
The ring $R$ acts on itself by multiplication.
For $f\in R$, let \defcolor{$m_f$} be the operator of multiplication by $f$: for $g\in R$, $m_f(g)=fg$.
By Stickelberger's Theorem~\cite{Cox2021}, the eigenvalues of $m_f$ are the values of $f$ at the  points of $\calV(I)$,
and the multiplicity of the eigenvalue $\lambda$ is the sum of the multiplicities in $\calV(I)$ of the inverse images,
$f^{-1}(\lambda)\cap\calV(I)$.

The (univariate) \defcolor{eliminant} \defcolor{$g$} of $I$ with respect to $f$ is the minimal polynomial of $m_f$.
When $f$ is a variable, e.g.\ $f=x_1$, $g$ is the monic generator of the univariate ideal $I\cap\KK[x_1]$.
In general, the eliminant generates the kernel of the map $\KK[Z]\to R$, where $Z\mapsto f$.
The function \texttt{regularRepresentation} computes a matrix representing $m_f$ with respect to the standard basis for $R$.
The function \texttt{univariateEliminant} returns the minimal polynomial of $m_f$, with respect to a new user-chosen variable
(or the default \texttt{Z}). 
\begin{leftbar}
\begin{verbatim}
i14 : S = QQ[x,y]

i15 : I = ideal(x^2*y^2-3*x^2-3*y^2+5,-3*x^2*y+2*x*y+4*x*y^2+1)

i16 : f = x + y 

i17 : regularRepresentation(f, I)

o17 = (| 1 x x2 xy xy2 y y2 y3 |, | 0 0 1/3   -1/3 -35/4 0 0 -105/16 |)
                                  | 1 0 5/3   0    0     0 0 -1/4    |
                                  | 0 1 -2/3  0    21/4  0 0 63/16   |
                                  | 0 1 -5/3  -2/3 0     1 0 -5/4    |
                                  | 0 0 19/9  7/3  1/2   0 1 23/24   |
                                  | 1 0 -20/9 0    1/4   0 0 -19/48  |
                                  | 0 0 -2/3  0    21/4  1 0 269/48  |
                                  | 0 0 4/3   0    0     0 1 1/2     |

o17 : Sequence

i18 : g = univariateEliminant(f, I)

           8       7         6        5          4         3          2
o18 = 1296Z  - 432Z  - 38223Z  - 4806Z  + 209784Z  + 14172Z  - 430242Z ...

o18 : QQ[Z]
\end{verbatim}
\end{leftbar}
The eliminant \defcolor{$g$} of $I$ with respect to $f$ defines the image of the scheme $\calV(I)$ under $f$.
When $g$ has degree equal to the degree of $I$, then $f$ is an isomorphism and thus $g$ may be used to study the scheme
$\calV(I)$.
For example, when both are reduced, $g$ and $R$ have the same Galois group over $\KK$.
While the eliminant {\it a priori} only tells us about $f(\calV(I))$, when $f$ is \demph{separating} (injective on the
points of $\calV(I)$), it tells us more about those points.
In our running example, both $g$ and $I$ have degree eight. (For $I$, note that this is the cardinality of the basis in \texttt{o17}.)
We see that $g$ is reduced and has four real roots.
\begin{leftbar}
\begin{verbatim}
i19 : T = ring(g)

i20 : gens gb ideal(g, diff(Z,g))

o20 = | 1 |

i21 : SturmCount(g)

o21 = 4
\end{verbatim}
\end{leftbar}
\noindent Thus $\calV(I)$ is reduced and consists of eight points, exactly four of which are real.

A useful variant of the eliminant is a \demph{rational univariate representation} of a zero-dimension\-al ideal $I$~\cite[Sect.\ 11.4]{BPR}.
This is a triple \defcolor{$(f,\chi,\phi)$} where $f$ is a linear form that is separating for $\calV(I)$, $\chi$ is the characteristic
polynomial of $m_f$---which retains the multiplicities of points of $\calV(I)$, if not their scheme structure---and $\phi$
is a rational map $\phi\colon\KK\to\KK^n$ that restricts to a bijection between the roots of $\chi$ and the points of $\calV(I)$.
\begin{leftbar}
\begin{verbatim}
i22 : (f, ch, ph) = rationalUnivariateRepresentation(I);

i23 : f

o23 = x + y

i24 : ch

       8   1 7   4247 6   89 5   8741 4   1181 3   71707 2   2051    6044
o24 = Z  - -Z  - ----Z  - --Z  + ----Z  + ----Z  - -----Z  - ----Z + ----
           3      144     24      54       108      216       324     27

i25 : ph

           7         6        5          4         3          2        ...
       864Z  + 21348Z  - 6066Z  - 231771Z  - 17610Z  + 701286Z  + 6986 ...
o25 = {--------------------------------------------------------------- ...
            7        6          5         4          3         2       ...
       5184Z  - 1512Z  - 114669Z  - 12015Z  + 419568Z  + 21258Z  - 430 ...

o25 : List
\end{verbatim}
\end{leftbar}
%

\section{Trace symmetric form}\label{S:three}
The remaining methods in \texttt{RealRoots} are linear-algebraic and may be used to count the points of $\calV(I)$ over any
field and to count real points of $\calV(I)$ according to 
the sign of another polynomial, similar to Sylvester's Theorem~\ref{Th:Sylvester}.
We demonstrate how this may be used for real root location.

A symmetric bilinear form \defcolor{$S$} in a real vector space $R$ has two basic invariants, its rank \defcolor{$\rho(S)$} and signature
\defcolor{$\sigma(S)$}.
If we choose a basis for $R$ and thus a corresponding matrix $M$ representing $S$, then $M$ will be symmetric, and therefore 
diagonalizable with all $d\vcentcolon = \dim(R)$ eigenvalues real.
The rank $\rho(M)$ of $M$ is its number of nonzero eigenvalues, and its signature  is the difference
\[
\sigma(M)\ \vcentcolon =\ \#\{\mbox{positive eigenvalues of $M$}\}
\ -\ \#\{\mbox{negative eigenvalues of $M$}\}\,.
\]
Sylvester's Law of Inertia asserts that the rank and signature are independent of the choice of basis, and therefore are invariants of the
symmetric form $S$.

Let $\KK$ be any field, let $I\subset\KK[x_1,\dotsc,x_n]$ be a zero-dimensional ideal, and set
$\defcolor{R}\vcentcolon=\KK[x_1,\dotsc,x_n]/I$, an  Artinian ring.    
For $f\in R$ (or in $\KK[x_1,\dotsc,x_n]$), multiplication by $f$ induces an endomorphism $m_f$ of $R$ as in Section~\ref{S:two}.
For $h\in R$ (or in $\KK[x_1,\dotsc,x_n]$), we define the symmetric bilinear \demph{trace form}, \defcolor{$S_h$} on $R$ by
$S_h(f,g)\vcentcolon= \mbox{trace}(m_{fgh})$.
The significance of the trace form is the following theorem.

\begin{theorem}[\protect{\cite[Thm.\ 4.72]{BPR}}]
  Suppose that $\KK$ is a subfield of $\RR$ and $I\subset\KK[x_1,\dotsc,x_n]$ is a zero-dimensional ideal with scheme
  $\calV(I)\subset\CC^n$.
  For $h\in\KK[x_1,\dotsc,x_n]$, the rank and signature of the trace form $S_h$ satisfy
  \begin{eqnarray}
    \rho(S_h)&=&    \label{Eq:traceCount}
      \#\{ z \in \calV(I) \mid h(z)\neq 0\} \quad and  \\
    \sigma(S_h)&=&
        \#\{ z \in \calV(I)\cap\RR \mid h(z)> 0\} \ -\  \nonumber
        \#\{ z \in \calV(I)\cap\RR \mid h(z)< 0\} \,.
  \end{eqnarray}
  If $\KK$ is any field with algebraic closure $\overline{\KK}$ and $\calV(I)$ is a subscheme of $\overline{\KK}^n$, then 
  the rank of the trace form $S_h$ satisfies~\eqref{Eq:traceCount}.
\end{theorem}

The rank of the trace form $S_1$ is used in \texttt{rationalUnivariateRepresentation} to certify that a linear form is separating.

\begin{figure}[htb]
  \centering
  \begin{picture}(168,120)(-3,0)
     \put(0,0){\includegraphics[height=120pt]{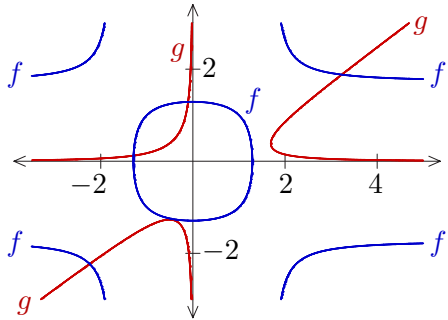}}
     \put( 2,4){{\color{myred}\small$g$}}
     \put(60,100){{\color{myred}\small$g$}}     \put(152,110){{\color{myred}\small$g$}}
     \put(-2,90){{\color{myblue}\small$f$}}     \put(158, 90){{\color{myblue}\small$f$}}
         \put(88,80){{\color{myblue}\small$f$}}
     \put(-2,26){{\color{myblue}\small$f$}}     \put(158,26){{\color{myblue}\small$f$}}
    \put(22,48){\small$-2$}    \put(101,48){\small$2$}   \put(135,48){\small$4$}
    \put(72,23){\small$-2$}    \put(72,92){\small$2$} 
  \end{picture}
  \caption{Two real curves}
  \label{F:two}
\end{figure}

We demonstrate how this may be used to study the number and location of real zeroes, using the ideal $I$ of Section~\ref{S:two}.
Let
\[
 {\color{myblue}f\ \vcentcolon=\ x^2y^2-3x^2-3y^2+5}
   \quad\mbox{and}\quad
  {\color{myred}g\ \vcentcolon=\ -3x^2y+4xy^2+2xy+1}\,.
\]
These define two curves in $\RR^2$, shown in Figure~\ref{F:two}.
While they have eight complex points in common, only four are real.
\begin{leftbar}
\begin{verbatim}
i26 : traceCount(I)

o26 = 8

i27 : realCount(I)

o27 = 4
\end{verbatim}
\end{leftbar}
\noindent We saw this in Section~\ref{S:two}, following \texttt{o21}.
When $h=1$, the rank and signature of the trace form $S_1$ count the complex and real points of $\calV(I)$, respectively,
and these are implemented as \texttt{traceCount(I)} and  \texttt{realCount(I)}.
Thus the possible tangency that we see in the third quadrant is only a near miss.

To see this in another way, consider the signature of the trace form $S_{y^2+2y}$,
\begin{leftbar}
\begin{verbatim}
i28 : signature traceForm(y^2 + 2*y, I)

o28 = 4
\end{verbatim}
\end{leftbar}
\noindent As this equals the number of real points of $\calV(I)$ and $y^2+2y<0$ for $-2<y<0$, there are no real points of
$\calV(I)$ in that horizontal strip, which includes the apparent `near tangency' in Figure~\ref{F:two}.

\end{document}